\documentclass[11pt]{article}
\usepackage[a4paper,margin=28mm]{geometry}
\usepackage{amsmath,amssymb,amsthm}
\usepackage[english]{babel}
\usepackage{hyperref}

\newtheorem{theorem}{Theorem}[section]
\newtheorem{lemma}[theorem]{Lemma}
\newtheorem{corollary}[theorem]{Corollary}
\newtheorem{result}[theorem]{Result}
\theoremstyle{remark}
\newtheorem{remark}[theorem]{Remark}

\newcommand{\Z}{\mathbb{Z}}
\newcommand{\Q}{\mathbb{Q}}
\newcommand{\F}{\mathbb{F}}
\newcommand{\Dih}{\mathrm{Dih}}
\DeclareMathOperator{\Aut}{Aut}

\title{Nonexistence of Hadamard difference sets in groups\\
 of order $4p^2$ with $p\equiv 3 \pmod 4$, $p\ge 7$}
\author{Bernhard Schmidt\\[1ex]
\small Division of Mathematical Sciences, School of Physical and Mathematical Sciences\\
\small Nanyang Technological University, 21 Nanyang Link, Singapore 637371\\
\small \texttt{bernhard@ntu.edu.sg}}
\date{}

\begin{document}
\maketitle

\begin{abstract}
We show that the Frobenius group $\langle x,y,z\mid x^p=y^p=z^4=1,\ xy=yx,\ zxz^{-1}=y,\ zyz^{-1}=x^{-1}\rangle$ of order $4p^2$ contains no $(4p^2,2p^2-p,p^2-p)$ difference set if $p\equiv 3\pmod 4$ is a prime with $p \geq 7$. Together with results of McFarland, Iiams and Liebler this shows that no group of order $4p^2$ with $p \equiv 3 \pmod 4$ and $p \geq 7$ contains a Hadamard difference set.
\end{abstract}


\section{Introduction}

Let $G$ be a finite group of order $v$. A $k$-subset $D$ of $G$ is called a $(v,k,\lambda)$ \emph{difference set} in $G$ if every nonidentity element of $G$ has exactly $\lambda$ representations $d_1d_2^{-1}$ with $d_1,d_2\in D$; the number $n=k-\lambda$ is the \emph{order} of $D$. Difference sets with parameters $(4u^2,2u^2-u,u^2-u)$ are called \emph{Hadamard difference sets}; they are equivalent to regular Hadamard matrices of order $4u^2$ admitting a regular group of automorphisms, see \cite{BJL,DJ}. We are concerned with the case $u=p$, $p$ an odd prime, i.e., with Hadamard difference sets in groups of order $4p^2$.

For $p=3$ nine of the fourteen groups of order $36$ contain $(36,15,6)$ difference sets \cite{Kib}. For $p\ge 5$ the Sylow $p$-subgroup $P$ of a group $G$ of order $4p^2$ is normal, so $G=P\rtimes K$ with $|K|=4$, and $G$ is determined up to isomorphism by $P$, $K$ and the conjugacy class of the image of $K$ in $\Aut(P)$. There are exactly $16$ such groups if $p\equiv 1\pmod 4$ and $12$ if $p\equiv 3\pmod 4$; we use the numbering $G_1,\dots,G_{16}$ of Iiams \cite{Iia}. By the results of McFarland \cite{McF}, Iiams \cite{Iia} and Liebler \cite{Lie}, the only group of order $4p^2$ with $p\equiv 3\pmod 4$, $p\ge 5$, for which the existence of a Hadamard difference set was open is the Frobenius group
\[
G_{13}=\langle x,y,z\mid x^p=y^p=z^4=1,\ xy=yx,\ zxz^{-1}=y,\ zyz^{-1}=x^{-1}\rangle
\]
of order $4p^2$, see \cite[Theorem 1.1]{Iia}. This case is settled by the following result.

\begin{theorem}\label{thm:G13}
Let $p\equiv 3\pmod 4$ be a prime with $p\ge 7$. Then $G_{13}$ contains no $(4p^2,2p^2-p,p^2-p)$ difference set.
\end{theorem}

\begin{corollary}\label{cor:main}
Let $p\equiv 3\pmod 4$ be a prime with $p\ge 7$. Then no group of order $4p^2$ contains a $(4p^2,2p^2-p,p^2-p)$ difference set.
\end{corollary}

To prove Theorem \ref{thm:G13}, we decompose a putative difference set $D$ into its intersections with the four cosets of $P$ in $G$ and use only the part of the group ring equation $DD^{(-1)}=n+\lambda G$ which lives in the subgroup $P\langle z^2\rangle$ of index $2$, the generalised dihedral group of $P$. This gives two equations in $\Z[P]$, namely a ``difference equation'' $\sum_j D_jD_j^{(-1)}=p^2+(p^2-p)P$ and a ``sum equation'' $D_0D_2+D_1D_3=\frac{p^2-p}{2}P$, in the spirit of the treatment of difference sets in dihedral groups by Leung, Ma and Wong \cite{LMW}. Applying a character of $P$ and comparing $(1-\zeta_p)$-adic valuations, these two equations lead to a system of two quadratic equations over $\F_p$ which has no nontrivial solution when $-1$ is a nonsquare modulo $p$. The argument is reminiscent of the proof of Eliahou, Kervaire and Saffari \cite{EKS} that the length of a Golay complementary pair, and hence of a Barker sequence of even length, has no prime factor $\equiv 3\pmod 4$.

\section{Preliminaries}\label{sec:prelim}

Let $G$ be a finite group. We identify a subset $X$ of $G$ with the element $\sum_{x\in X}x$ of the integral group ring $\Z[G]$. For $A=\sum_{g\in G}a_gg\in\Z[G]$ we write $A^{(-1)}=\sum_ga_gg^{-1}$ and $|A|=\sum_ga_g$. A $k$-subset $D$ of $G$ is a $(v,k,\lambda)$ difference set in $G$ if and only if
\begin{equation}\label{eq:ds}
DD^{(-1)}=n+\lambda G\quad\text{in }\Z[G],
\end{equation}
where $n=k-\lambda$. Characters of abelian groups are extended to the group rings by linearity. If $\chi$ is a nontrivial linear character of $G$, then applying $\chi$ to \eqref{eq:ds} gives $|\chi(D)|^2=n$, since $\chi(G)=0$.

For a positive integer $m$ let $\zeta_m$ be a primitive complex $m$th root of unity. We need the following two facts.

\begin{result}\label{res:kronecker}
If $x_1,\dots,x_r$ are totally nonnegative algebraic integers with $x_1+\dots+x_r=1$, then exactly one of the $x_i$ is $1$ and the others are $0$.
\end{result}
\begin{proof}
All conjugates of each $x_i$ lie in $[0,1]$, so the norm of $x_i$ is an integer in $[0,1]$. If the norm is $0$, then $x_i=0$; if it is $1$, then all conjugates of $x_i$ are $1$, so $x_i=1$.
\end{proof}

\begin{result}[Ma's lemma \cite{Ma}]\label{res:ma}
Let $p$ be a prime, let $H$ be a finite abelian group with cyclic Sylow $p$-subgroup, and let $U$ be the unique subgroup of order $p$ of $H$. If $Y\in\Z[H]$ satisfies $\chi(Y)\equiv 0\pmod p$ for all nontrivial characters $\chi$ of $H$, then $Y=pX_1+UX_2$ for some $X_1,X_2\in\Z[H]$.
\end{result}

Let $p$ be an odd prime and $\zeta=\zeta_p$. Then $\mathfrak p=(1-\zeta)$ is the unique prime ideal of $\Z[\zeta]$ above $p$, and $(p)=\mathfrak p^{p-1}$. Let $\nu$ denote the $\mathfrak p$-adic valuation on $\Z[\zeta]$, so $\nu(p)=p-1$ and $\nu\bigl((1-\zeta)(1-\bar\zeta)\bigr)=2$. The residue field $\Z[\zeta]/\mathfrak p$ is $\F_p$. Every Galois automorphism $\sigma$ of $\Q(\zeta)$ fixes $\mathfrak p$ and acts trivially on $\Z[\zeta]/\mathfrak p$, since $\sigma(\zeta)\equiv 1\equiv\zeta\pmod{\mathfrak p}$; in particular $\bar x\equiv x\pmod{\mathfrak p}$ and $\nu(\bar x)=\nu(x)$ for all $x\in\Z[\zeta]$.

From now on let $p\ge 5$, $G=G_{13}$, $P=\langle x,y\rangle$, and let $D$ be a $(4p^2,2p^2-p,p^2-p)$ difference set in $G$; thus $n=p^2$ and $\lambda=p^2-p$. Every $A\in\Z[G]$ can be written uniquely as $A=A_0+A_1z+A_2z^2+A_3z^3$ with $A_j\in\Z[P]$; we call $A_jz^j$ the \emph{component} of $A$ in $\Z[P]z^j$. In particular,
\begin{equation}\label{eq:blocks}
D=D_0+D_1z+D_2z^2+D_3z^3,\qquad D_j=(D\cap Pz^j)z^{-j}\subseteq P ,
\end{equation}
and we call the $D_j$ the \emph{blocks} of $D$. Since $(D_jz^j)(D_kz^k)^{(-1)}=D_j\bigl(z^{j-k}D_k^{(-1)}z^{k-j}\bigr)z^{j-k}$ lies in $\Z[P]z^{j-k}$, the component of $DD^{(-1)}$ in $\Z[P]$ is $\sum_jD_jD_j^{(-1)}$, and comparing the components in $\Z[P]$ of both sides of \eqref{eq:ds} gives
\begin{equation}\label{eq:F}
D_0D_0^{(-1)}+D_1D_1^{(-1)}+D_2D_2^{(-1)}+D_3D_3^{(-1)}=p^2+(p^2-p)P .
\end{equation}

By a \emph{line} of $P$ we mean a subgroup $L$ of $P$ of order $p$. For $A\in\Z[P]$ let $A_L\in\Z[P/L]$ denote the image of $A$ under the natural epimorphism $\Z[P]\to\Z[P/L]$; the coefficient of $A_L$ at a coset $xL$ is the sum of the coefficients of $A$ over $xL$, so $A_L$ has coefficients in $\{0,\dots,p\}$ if $A$ is a subset of $P$. The characters of $P$ with kernel $L$ are the nontrivial characters of $P/L\cong C_p$; they are Galois conjugate to each other, and $\chi(A)=\chi(A_L)$ for these characters.

\begin{lemma}\label{lem:sizes}
If $p\equiv 3\pmod 4$, then the multiset of block sizes $\{|D_0|,|D_1|,|D_2|,|D_3|\}$ is $\{p(p+1)/2,\,p(p-1)/2,\,p(p-1)/2,\,p(p-1)/2\}$.
\end{lemma}

\noindent Lemma \ref{lem:sizes} is \cite[Lemma 5.3]{Iia}; it follows from the self-conjugacy of $p$ modulo $4$ applied to the characters of $G/P\cong C_4$.

\begin{lemma}\label{lem:spikes}
Let $p\equiv 3\pmod 4$, $p\ge 7$, and let $\chi$ be a nontrivial character of $P$. Then $p$ does not divide $\chi(D_j)$ in $\Z[\zeta]$ for at least one $j$.
\end{lemma}
\begin{proof}
Suppose that $p$ divides $\chi(D_j)$ for all $j$, and let $L=\ker\chi$. Since the nontrivial characters of $P/L$ are Galois conjugate to $\chi$, $p$ divides $\psi((D_j)_L)$ for all nontrivial characters $\psi$ of $P/L$ and all $j$. By Result \ref{res:ma} applied to $H=P/L$, we have $(D_j)_L=pX_j+c_j(P/L)$ with $X_j\in\Z[P/L]$ and $c_j\in\Z$, and we may assume $0\le c_j<p$. As the coefficients of $(D_j)_L$ lie in $\{0,\dots,p\}$, we conclude that $X_j=0$ if $c_j>0$, and that $X_j$ is a subset of $P/L$ if $c_j=0$. In both cases $\chi(D_j)=p\chi(X_j)$. Applying $\chi$ to \eqref{eq:F} gives $\sum_j|\chi(D_j)|^2=p^2$, hence $\sum_j|\chi(X_j)|^2=1$. By Result \ref{res:kronecker}, there is $j_0$ with $|\chi(X)|^2=1$ for $X=X_{j_0}$. Then $X\ne 0$, so $c_{j_0}=0$, $X$ is a subset of $P/L$ and $|D_{j_0}|=p|X|$. Moreover, $|\psi(X)|^2=1$ for all nontrivial characters $\psi$ of $P/L$, so $XX^{(-1)}-1$ is annihilated by all nontrivial characters of $P/L$ and hence $XX^{(-1)}=1+c\,(P/L)$ for some integer $c$. Comparing the coefficients of the identity element and the sums of all coefficients on both sides gives $|X|=1+c$ and $|X|^2=1+cp$, so $(|X|-1)(|X|-p+1)=0$ and $|D_{j_0}|\in\{p,p^2-p\}$. By Lemma \ref{lem:sizes}, $|D_{j_0}|\in\{p(p+1)/2,p(p-1)/2\}$, which forces $p\le 3$.
\end{proof}

\section{Proof of Theorem \ref{thm:G13}}\label{sec:proof}

Let $p\equiv 3\pmod 4$, $p\ge 7$, and let $D$ be a $(4p^2,2p^2-p,p^2-p)$ difference set in $G=G_{13}$, written as in \eqref{eq:blocks}. We only use that $z^2$ acts on $P$ by inversion, i.e., $z^2Az^{-2}=A^{(-1)}$ for all $A\in\Z[P]$; indeed, $z^2xz^{-2}=zyz^{-1}=x^{-1}$ and $z^2yz^{-2}=zx^{-1}z^{-1}=y^{-1}$.

The component of $DD^{(-1)}$ in $\Z[P]z^2$ is the sum of the terms $D_j\bigl(z^{j-k}D_k^{(-1)}z^{k-j}\bigr)z^{j-k}$ with $j-k\equiv 2\pmod 4$. As $z^{j-k}=z^2$ and $z^2D_k^{(-1)}z^{-2}=D_k$, each of these terms equals $D_jD_kz^2$, and comparing the components in $\Z[P]z^2$ of both sides of \eqref{eq:ds} gives $2(D_0D_2+D_1D_3)z^2=(p^2-p)Pz^2$, i.e.,
\begin{equation}\label{eq:S}
D_0D_2+D_1D_3=\tfrac{p^2-p}{2}\,P .
\end{equation}
Let $\chi$ be a nontrivial character of $P$ and $\beta_j=\chi(D_j)\in\Z[\zeta]$. Applying $\chi$ to \eqref{eq:F} and \eqref{eq:S} gives
\begin{equation}\label{eq:FS}
\sum_{j=0}^3|\beta_j|^2=p^2\qquad\text{and}\qquad\beta_0\beta_2+\beta_1\beta_3=0 .
\end{equation}
The contradiction will come from the following observation.

\begin{lemma}\label{lem:Fp}
Let $p\equiv 3\pmod 4$. There are no $u_0,u_1,u_2,u_3\in\F_p^*$ with
\[
u_0^2+u_1^2+u_2^2+u_3^2=0\qquad\text{and}\qquad u_0u_2+u_1u_3=0 .
\]
\end{lemma}
\begin{proof}
The two equations give $(u_0+u_2)^2+(u_1+u_3)^2=0$ and $(u_0-u_2)^2+(u_1-u_3)^2=0$. Since $-1$ is a nonsquare in $\F_p$, a sum of two squares in $\F_p$ vanishes only if both squares vanish. Hence $u_2=-u_0$ and $u_2=u_0$, so $2u_0=0$, a contradiction.
\end{proof}

\begin{proof}[Proof of Theorem \ref{thm:G13}]
Fix a nontrivial character $\chi$ of $P$ and let $\beta_j=\chi(D_j)$ as above. By \eqref{eq:FS}, not all $\beta_j$ vanish. For $\beta_j\ne 0$ put $\nu_j=\nu(\beta_j)$ and note that $\nu(|\beta_j|^2)=2\nu_j$. Let $m$ be the minimum of the $\nu_j$ over all $j$ with $\beta_j\ne 0$, and let $I=\{j:\beta_j\ne 0,\ \nu_j=m\}$. We claim that
\begin{equation}\label{eq:claim}
m<p-1\qquad\text{and}\qquad I=\{0,1,2,3\} .
\end{equation}
If $m\ge p-1$, then $p$ divides all $\beta_j$, contradicting Lemma \ref{lem:spikes}. If $|I|=1$, then $\sum_j|\beta_j|^2$ has valuation $2m<2(p-1)=\nu(p^2)$, contradicting \eqref{eq:FS}. Now write $\beta_j=(1-\zeta)^m\beta_j'$ with $\beta_j'\in\Z[\zeta]$, so that $\beta_j'$ is a unit for $j\in I$ and $\beta_j'\in\mathfrak p$ for $j\notin I$, and let $u_j\in\F_p$ be the residue of $\beta_j'$ modulo $\mathfrak p$. Dividing the first equation of \eqref{eq:FS} by $\bigl((1-\zeta)(1-\bar\zeta)\bigr)^m$ gives
\begin{equation}\label{eq:div}
\sum_{j=0}^3\beta_j'\bar\beta_j'=\frac{p^2}{\bigl((1-\zeta)(1-\bar\zeta)\bigr)^m},
\end{equation}
where the right hand side is an element of $\Z[\zeta]$ of valuation $2(p-1)-2m>0$. Reducing \eqref{eq:div} modulo $\mathfrak p$ and using $\bar\beta_j'\equiv\beta_j'$ we obtain
\begin{equation}\label{eq:sumsq}
\sum_{j\in I}u_j^2=0\quad\text{in }\F_p,\qquad u_j\ne 0\text{ for }j\in I .
\end{equation}
If $|I|=2$, then \eqref{eq:sumsq} is a vanishing sum of two nonzero squares in $\F_p$, which is impossible for $p\equiv 3\pmod 4$. If some $\beta_j$ vanishes, then the second equation of \eqref{eq:FS} shows that a second $\beta_{j'}$ vanishes, so $|I|\le 2$. Hence, if $|I|=3$, all $\beta_j$ are nonzero, and taking valuations in the second equation of \eqref{eq:FS} gives $\nu_0+\nu_2=\nu_1+\nu_3$; but one side equals $2m$ and the other is larger than $2m$. This proves \eqref{eq:claim}.

By \eqref{eq:claim}, all $\beta_j'$ are units, and \eqref{eq:sumsq} reads $\sum_{j=0}^3u_j^2=0$ with all $u_j\ne 0$. Moreover, dividing the second equation of \eqref{eq:FS} by $(1-\zeta)^{2m}$ gives $\beta_0'\beta_2'+\beta_1'\beta_3'=0$, hence $u_0u_2+u_1u_3=0$. This contradicts Lemma \ref{lem:Fp}.
\end{proof}

\begin{proof}[Proof of Corollary \ref{cor:main}]
As stated in the introduction, $G_{13}$ is the only group of order $4p^2$ with $p\equiv 3\pmod 4$, $p\ge 5$, which is not covered by the results of McFarland \cite{McF}, Iiams \cite{Iia} and Liebler \cite{Lie}, and $G_{13}$ is excluded by Theorem \ref{thm:G13}.
\end{proof}

\section{Remarks}

\begin{remark}\label{rem:dihedral}
Let $N=P\langle z^2\rangle$, $X=D\cap N$ and $Y=(D\cap Nz)z^{-1}$. Adding the components of \eqref{eq:ds} in $\Z[P]$ and in $\Z[P]z^2$ gives
\begin{equation}\label{eq:index2}
XX^{(-1)}+YY^{(-1)}=p^2+(p^2-p)N\quad\text{in }\Z[N],
\end{equation}
and \eqref{eq:F}, \eqref{eq:S} are the two components of this equation; $N$ is the generalised dihedral group $\Dih(P)$ of $P=C_p^2$. The proof of Theorem \ref{thm:G13} uses only \eqref{eq:index2} and the block sizes, and the latter follow from \eqref{eq:index2} via the sign character of $N/P$ once $\{|X|,|Y|\}=\{p^2,p^2-p\}$ is known. Hence the argument shows: \emph{for $p\equiv 3\pmod 4$, $p\ge 7$, there are no subsets $X,Y$ of $\Dih(C_p^2)$ with $\{|X|,|Y|\}=\{p^2,p^2-p\}$ and $XX^{(-1)}+YY^{(-1)}=p^2+(p^2-p)\Dih(C_p^2)$.} Consequently, every group of order $4p^2$ containing $\Dih(C_p^2)$ as a subgroup of index $2$ has no Hadamard difference sets for $p\equiv 3\pmod 4$, $p\ge 7$; besides $G_{13}$ these are the groups $G_9=\Dih(C_p^2)\times C_2$ and $G_{10}=D_{2p}\times D_{2p}$, which are covered by \cite[Theorem 7.1]{Iia}. For a single subset this is the setting of \cite{LMW}.
\end{remark}

\begin{remark}
In the extended version \cite{Sch} of this paper, the block decomposition \eqref{eq:blocks} is also used to give self-contained proofs of Iiams' results \cite[Theorems 7.1 and 7.2]{Iia} for all primes $p\ge 5$ and of the results of McFarland \cite{McF}, Iiams \cite{Iia} and Liebler \cite{Lie} for $p\equiv 3\pmod 4$; thus Corollary \ref{cor:main} is proved there without using any of these results. In particular, for $p\equiv 3\pmod 4$ the use of Segre's theorem on arcs in $PG(2,p)$ in Liebler's proof for the group $C_p^2\rtimes C_4$, where $C_4$ acts by inversion, can be avoided.
\end{remark}

\begin{remark}
The condition $p\ge 7$ is used only in Lemma \ref{lem:spikes}: for $p=3$ the block sizes $p=3$ and $p^2-p=6$ obtained in its proof are admissible. Indeed, the group $G_{13}$ of order $36$ contains the $(36,15,6)$ difference set $D=D_0+D_1z+D_2z^2+D_3z^3$ with
\[
D_0=\langle x\rangle+y\langle x\rangle,\quad D_1=\{x,\,x^2y,\,y^2\},\quad D_2=\{x,\,xy,\,xy^2\},\quad D_3=\{x,\,y,\,x^2y^2\}.
\]
Here $\chi(D_1)=\chi(D_2)=\chi(D_3)=0$ and $\chi(D_0)=3(1+\chi(y))$ for the characters $\chi$ of $P$ with kernel $\langle x\rangle$, so $3$ divides all $\chi(D_j)$. 
\end{remark}

\begin{remark}
For $p\equiv 1\pmod 4$, the existence of Hadamard difference sets is undecided only for the groups $G_4$, $G_{13}$, $G_{15}$ and $G_{16}$: all other groups of order $4p^2$ are excluded by \cite[Theorem 1.1]{Iia} and by AbuGhneim \cite{Abu}, who settled $G_{11}$ and $G_{14}$. Equations \eqref{eq:FS} hold in every group $C_p^2\rtimes\langle z\rangle$ in which $z^2$ acts by inversion, in particular in $G_{13}$ and $G_{16}=C_p^2\rtimes C_4$ ($z$ acting as a scalar of order $4$) for $p\equiv 1\pmod 4$. For such $p$, however, $-1$ is a square in $\F_p$, say $i^2=-1$, and the system in Lemma \ref{lem:Fp} has nontrivial solutions such as $(u_0,u_1,u_2,u_3)=(1,i,-1,-i)$. This is in accordance with the existence of $(100,45,20)$ difference sets in nonabelian groups of order $100$ \cite{Smi,GV}.
\end{remark}

\section*{Declaration on the use of generative AI}
The research reported in this paper was carried out under the direction of the author with the assistance of the AI model Claude (Anthropic, Fable 5.1). The author posed the problem, chose the method, namely the decomposition of a putative difference set into the blocks \eqref{eq:blocks} and the resulting equations \eqref{eq:F} and \eqref{eq:S}, in analogy to the treatment of difference sets in dihedral groups \cite{LMW}, and directed the investigation. The AI model found the proof of Theorem \ref{thm:G13} including Lemma \ref{lem:Fp}, and produced a first draft of the manuscript. The author has checked all arguments in detail and takes full responsibility for the content of the paper.

\end{document}